\documentstyle[12pt]{article}
\newcommand{\fr}{\frac}
\newcommand{\lb}{\label}
\newcommand{\ti}{\tilde}
\newcommand{\be}{\begin{equation}}
\newcommand{\ee}{\end{equation}}
\newcommand{\ba}{\begin{array}}
\newcommand{\ea}{\end{array}}
\newcommand{\beqa}{\begin{eqnarray}}
\newcommand{\al}{\alpha}

\newcommand{\del}{\partial}
\newcommand{\eeqa}{\end{eqnarray}}

\newcommand{\Pc}{{\cal P}}
\newcommand{\Oc}{{\cal O}}
\newcommand{\Cc}{{\cal C}}

\newcommand{\Qc}{{\cal Q}}
\newcommand{\Sc}{{\cal S}}
\newcommand{\Sh}{{\cal S}^{(\hbar )}}
\newcommand{\Sq}{{\cal S}^{(q)}}

\newcommand{\th}{T^{(\hbar )}}
\newcommand{\tq}{T^{(q )}}
\newcommand{\tbir}{T^{(GF )}}

\newcommand{\ph}{P_{\hbar}}
\newcommand{\xh}{X_{\hbar}}
\newcommand{\pq}{P_{q}}
\newcommand{\xq}{X_{q}}

\newcommand{\dl}{\stackrel{\leftarrow}{\partial}}
\newcommand{\dr}{\stackrel{\rightarrow}{\partial}}
\newcommand{\Dl}{\stackrel{\leftarrow}{D}}
\newcommand{\Dr}{\stackrel{\rightarrow}{D}}

\begin{document}

\begin{flushright}
MRC.PH.TH--12/96

q-alg/9609023
\end{flushright}

\vspace{1cm}

{\Large {\bf 
\noindent
q--DEFORMED STAR PRODUCTS AND MOYAL BRACKETS } }

\vspace{1.5cm}

\noindent
\"{O}mer F. DAYI\footnote{E-mail address:
dayi@mam.gov.tr.}, 

\vspace{1cm}                                                                   

{\small {\it
\noindent
T\"{U}B\.{I}TAK-Marmara Research Centre,
Research Institute for Basic Sciences,
\mbox{Department} of Physics,
P.O.Box 21, 41470 Gebze, Turkey.  } }

\vspace{2cm}

{\small
\begin{center}
{\bf Abstract}
\end{center}

\vspace{.2cm}

The standard and anti--standard
ordered operators
acting on two--dimensional q--deformed phase space
are shown to satisfy algebras 
which can be called \mbox{q--$W_\infty .$}
q--star products and q--Moyal
brackets  corresponding to these
algebras are constructed. Some applications like
defining q--classical mechanics
and q--path integrals are discussed.

\vspace{2cm}

PACS: 03.65 Fd; 03.65 Ca.

\vspace{1cm}

\pagebreak

\section{Introduction }

A way of visualizing quantum mechanics as $\hbar$--deformation of the
classical case is to utilize symbols of operators, star
products and Moyal brackets\cite{bea}.
If we denote the symbol map
by $\Sh_O$, to an operator $\hat{f}$ there corresponds a
c--number object $f^{(\hbar )}_O:$
\[
\Sh_O\left(\hat{f}\right)= f^{(\hbar )}_O,
\]
where the subscript $O$ denotes the operator ordering
adopted. Non-commutativity of quantum mechanics is
taken into account in terms of the star product defined
such that
\[
\Sh_O\left(\hat{f}\hat{g}\right)=
\Sh_O\left(\hat{f}\right)\star^{\hbar}_O 
\Sh_O\left(\hat{g}\right)
\]
and the symbol corresponding to the commutator
(divided by $\hbar$)
of the operators $\hat{f}$ and $\hat{g}$
is the Moyal bracket (or $\star^\hbar$--bracket):
\be
\lb{moy}
\{ f^{(\hbar )}_O(p,x), g^{(\hbar )}_O(p,x)\}_M
\equiv \fr{1}{i\hbar}
\left( f^{(\hbar )}_O(p,x)\star^\hbar_O g^{(\hbar )}_O(p,x)
- g^{(\hbar )}_O(p,x)
\star^\hbar_O f^{(\hbar )}_O(p,x) \right),
\ee
where $p$ and $x$  are classical  phase space variables.

Recently $\star$--products and Moyal brackets are used in 
some diverse context like  
constructing a representation of non--commuting forms\cite{ma}
and they appeared in the formulation
of two--dimensional area preserving diffeomorphisms\cite{apd}.
The latter is closely related to the quantum Hall
effect which also exhibits a q--deformed structure
where q is related to filling fraction\cite{qhe}.
Hence, understanding q--deformations in terms of
symbols of the q--deformed operators
may shed light on some aspects of the quantum Hall effect.
Although the
quantum Hall effect is a many body problem,
it is based on  a field theory which is given 
in a two--dimensional phase space when it is restricted to the
lowest Landau level.
Hence, understanding q--deformed star products and
Moyal brackets in a two--dimensional phase space
is the first step in this direction.

One of the essential properties of the ordinary star 
product is its associativity. 
In this case it is known that
the unique deformation of the Poisson bracket is the Moyal
bracket  (\ref{moy})\cite{umb}. On the other hand it was
shown that it is impossible to deforme Heisenberg dynamics 
by keeping the algebra of the observables 
associative\cite{pil}.
Hence, to obtain a nontrivial deformation of the 
quantum mechanics we should sacrifice associativity
of the star product. In fact, 
q--deformed algebras are not associative, so that 
a q--deformed star product leading to them should also
be non-associative.

Here,  after deriving the algebras of
the standard (or $XP$) and the anti--standard 
(or $PX$) ordered operators
q--deformed star products and Moyal brackets
corresponding to them
are constructed and the obstacles 
in generalizing Weyl ordering to 
the q--deformed operators are emphasized.

Once a q--Moyal bracket is obtained
an immediate application  is to define
q--classical mechanics in terms of q--Poisson brackets.
This and some other applications like classical
as well as quantum mechanical properties of q--canonical
transformations and a definition of q--path 
integrals on general grounds are  briefly discussed.

\section{Ordinary Symbols and Star Products:}

In the two--dimensional phase space given in terms of 
the ordinary ($\hbar$--deformed)
canonical operators 
\be
[\ph , \xh ] =i\hbar ,
\ee
an operator $\Oc$ can be written as
\be 
\lb{gom}
\Oc =\sum_{(m,n)>0} \Oc_{m,n} g(\ph^m,\xh^n),
\ee
where $g(\ph^m,\xh^n)$ are some functions of $\ph^m,\xh^n$
and $\Oc_{m,n}$ are some constant coefficients
depending on the operator ordering scheme adopted. Hence, it
is sufficient to deal with the monomials $g(\ph^m,\xh^n)$
as far as the algebraic properties of the operators
are concerned. 
We are interested in the following  
three different ordering schemes which constitute  complete
basis.

Standard (or $XP $) ordering: In this scheme one deals with
the monomials $\xh^n \ph^m$. They  satisfy the 
Lie algebra\cite{du}
\be
\lb{sal}
[\xh^n \ph^m ,\xh^k \ph^l ]= 
\sum_{r=1}^\infty (i\hbar )^r r! \left\{
\left(
\ba{c}
k \\
r
\ea
\right)
\left(
\ba{c}
m \\
r
\ea
\right)
-\left(
\ba{c}
n \\
r
\ea
\right)
\left(
\ba{c}
l \\
r
\ea
\right) \right\}
\xh^{n+k-r} \ph^{l+m-r}.
\ee

Antistandard (or $ PX$) ordering:  the suitable
monomials $\ph^n \xh^m,$ can be shown to satisfy\cite{du}
\be
\lb{aal}
[\ph^n \xh^m ,\ph^k \xh^l ]= 
\sum_{r=1}^\infty (-i\hbar )^r r!\left\{
\left(
\ba{c}
k \\
r
\ea
\right)
\left(
\ba{c}
m \\
r
\ea
\right)
-\left(
\ba{c}
n \\
r
\ea
\right)
\left(
\ba{c}
l \\
r
\ea
\right) \right\}
\ph^{n+k-r} \xh^{l+m-r}.
\ee

The algebras (\ref{sal}) and (\ref{aal}) 
are called $W_\infty $\cite{apd}.

Weyl ordering: The monomials 
\be
\lb{hsm}
\th_{m,n} =\exp (\ph \fr{\del}{\del \Pc} 
+ \xh \fr{\del}{\del \Qc}) \Pc^m\Qc^n|_{\Pc=\Qc=0},
\ee
where $\Pc$ and $\Qc$ are c--number variables, constitute
a complete basis and satisfy the Lie algebra\cite{bd},\cite{gf}
\be
\lb{la1}
[ \th_{m,n}, \th_{k,l} ] =\sum_{a=0}^b
(i\hbar )^{2a+1} B_{mnkl}^a \th_{m+k-2a-1,n+l-2a-1},
\ee
where $b=\min \{ [(m+k-1)/2 ,(n+l-1)/2]\}.$
$B_{00kl}=B_{mn00}\equiv 0$ and for the other values of the 
indices
\be
B_{mnkl}^a\equiv \sum_{c=0}^{2a+1}\fr{(-1)^c m! n! k! l!}{
(2a+1-b)! b! (m+c-2a-1)! (n-c)! (k-c)! (l+c-2a-1)! } .
\ee

For our purposes it is sufficient to deal with the
symbols of the monomials in each ordering scheme,
which are given by
\be
\Sh_S[\xh^m \ph^n]=
\Sh_A[\ph^m \xh^n]=
\Sh_W[ \th_{m,n}]= p^mx^n,
\ee
where $p$ and $x$ are c--number variables.

Related star products are\cite{bea},\cite{du}
\beqa
\star_S^\hbar & \equiv & \exp \left( i\hbar 
\dl_p \dr_x
\right) , \lb{ssp} \\
\star_A^\hbar & \equiv & \exp \left( -i\hbar 
\dl_x \dr_p
\right) , \lb{asp} \\
\star_W^\hbar & \equiv &
\exp \left[ \fr{-i\hbar}{2}
(\dl_x \dr_p
-\dl_p\dr_x)\right] .          \lb{wsp}
\eeqa
We used $\fr{\del}{\del x}\equiv \del_x;\ 
\fr{\del}{\del p}\equiv \del_p.$
One of them can be utilized in the Moyal bracket (\ref{moy})
in terms of the related symbols.
Observe that $\star^\hbar$--products are associative\cite{bea},
so that the algebraic properties of the commutators
are preserved.

\section{Classical $(\hbar =0)$ q--Deformed Symbols, 
Star Products and Moyal Brackets:}

Classical $(\hbar =0)$ q--deformed canonical operators 
are defined as
\be
\pq \xq-q\xq \pq =0.
\ee

Now, the monomials of the different ordering schemes are
equivalent up to an overall q--dependent factor. Thus
it is sufficient to consider the algebra 
\be
\lb{alq}
q^{nk}\tq_{m,n}\tq_{k,l}
-q^{ml}\tq_{k,l}\tq_{m,n}=0,
\ee                         
where  $\tq_{m,n}\equiv \pq^m \xq^n.$ 
They give a complete basis in $\pq ,\xq$ space and 
the basis operators in another ordering scheme 
are equal to $\tq_{m,n}$ up to an overall q--dependent constant.

Symbols of the operators $\tq_{m,n}$ are
\be
\lb{s2}
\Sq \tq_{m,n}= p^mx^n.
\ee

It is possible to define  associative 
$\star^q$--products
\beqa
\star^q_S & \equiv &
\exp \left(\nu \dl_p p x\dr_x\right) , \lb{nus} \\
\star^q_A & \equiv &
\exp \left(-\nu 
\dl_xxp\dr_p\right), \\
\star^q_W & \equiv &
\exp \left[ \fr{-\nu }{2}
(\dl_xxp\dr_p
-\dl_ppx\dr_x)\right] , \lb{wus}
\eeqa
where $\nu\equiv \ln q$ 
and the subscripts denote the resemblance to the ordinary
star products (\ref{ssp})--(\ref{wsp}).
In fact, (\ref{nus})--(\ref{wus}) can be obtained from
(\ref{ssp})--(\ref{wsp}) by a canonical transformation
in accord with the fact that there exists a unique deformation
of the Poisson bracket if the star product is assocaitive. 

The symbols defined in 
(\ref{s2}) satisfy 
\be
q^{nk}p^mx^n\star^q p^kx^l
-q^{ml}p^kx^l\star^q p^mx^n=0,
\ee
for any $\star^q$--product  (\ref{nus})--(\ref{wus}).

\section{q--symbols, q--star Products and q--Moyal Brackets: }

q--deformed $(\hbar \neq 0)$ canonical variables
satisfy 
\be
PX-qXP=i\hbar.
\ee
To reproduce (\ref{alq}) in $\hbar =0$ limit we 
define the q--commutator as
\be
\lb{qc}
[t(P^m,X^n),t(P^k,X^l)]_q\equiv
q^{nk} t(P^m,X^n)t(P^k,X^l)-
q^{ml}t(P^k,X^l)t(P^m,X^n),
\ee
where $t(P^m,X^n)$ is a function of $P^m$ and
$X^n,$ depending on the operator ordering scheme adopted.
Observe that the weights of the 
q--commutator (\ref{qc}) change according 
to the operators which are considered. 

As one can easily see,
the standard or the anti--standard ordered monomials
form a complete basis in the q--phase space
given by $P$ and $X.$

By explicit calculations 
q--algebra satisfied by
the standard ordered monomials $X^nP^m$ 
can be derived:
\beqa
[X^n P^m ,X^k P^l ]_q & = & 
\sum_{r=1}^\infty (i\hbar )^r [r]! \left\{ 
q^{(k-r)(n-r)+ml}
\left[
{\ba{c}
k \\
r
\ea  }
\right]
\left[
{\ba{c}
m \\
r
\ea}
\right]  
 -q^{(m-r)(l-r)+nk}
\left[
\ba{c}
n \\
r
\ea
\right]
\left[
\ba{c}
l \\
r
\ea
\right] \right\} \nonumber \\
& & X^{n+k-r} P^{l+m-r}, \lb{qsal}
\eeqa
where we used the definitions of the q--factorial
\[
[n]!\equiv \left(\fr{1-q^n}{1-q}\right) !
\equiv [1] [2]\cdots [n-1][n],
\]
and the q--binomial coefficient
\[
\left[
\ba{c}
n \\
r
\ea
\right] \equiv
\fr{[n]!}{[n-r]![r]!} .
\]

Similarly we see that
in the anti--standard ordering scheme the
monomials $P^mX^n$ satisfy the q--algebra
\be
\lb{qaal}
[P^n X^m ,P^k X^l ]_q= 
\sum_{r=1}^\infty (-i\hbar )^r q^{r(r-1)/2}  [r]!\left\{
\left[
\ba{c}
k \\
r
\ea
\right]
\left[
\ba{c}
m \\
r
\ea
\right]
-\left[
\ba{c}
n \\
r
\ea
\right]
\left[
\ba{c}
l \\
r
\ea
\right] \right\}
P^{n+k-r} X^{m+l-r}.
\ee

The q--deformed algebras (\ref{qsal}) and (\ref{qaal}) 
can be called q--$W_\infty .$
For other definitions of q--$W_\infty$ algebra see \cite{qwi}.

Symbol maps for the standard and anti-standard
orderings are defined like the ordinary case:
\be
\lb{symq}
\Sc_S(X^m P^n)=
\Sc_A(P^m X^n)=
p^mx^n.
\ee

As emphasized before,
the ordinary star products should be
associative, so that, the Moyal bracket satisfies
an identity corresponding to the Jacobi identity satisfied by
the ordinary operators.  However, now the situation is
altered drastically: we do not any more deal with the commutators
which are not aware of their entries but with the q--commutators
(\ref{qc}) which
change according to their entries i.e. the underlyning algebraic 
structure is non-associative. 
Thus the associativity
condition cannot be preserved and we can obtain a non-trivial 
deformation of the Poisson bracket other than $\hbar$--deformation.

In terms of the q--derivative
\be
D_zf(z)\equiv \fr{f(z)-f(qz)}{(1-q)z},
\ee
we can construct q--star products for the standard and anti--standard
orderings as
\beqa
\star_S & \equiv & \sum_{r=0}^\infty \fr{(i\hbar )^r}{[r]!}
\Dl_p^r \exp (\nu \dl_ppx\dr_x) \Dr_x^r, \lb{qsso} \\
\star_A & \equiv & \sum_{s=0}^\infty (-\nu \dl_xx)^s
\sum_{r=0}^\infty \fr{(-i\hbar )^r q^{r(r-1)/2}}{[r]!} 
\Dl_x^r  \Dr_p^r (p\dr_p)^s.             \lb{qaso} 
\eeqa

One can see that if the q--Moyal bracket is defined as
\be
\lb{qmoy}
\{p^mx^n,p^kx^l\}_{q-M}\equiv
\fr{1}{i\hbar }(q^{nk} p^mx^n\star p^kx^l-
q^{ml}p^kx^l\star p^mx^n),
\ee
by using  the symbols  (\ref{symq}) and  
the q--star products
(\ref{qsso})--(\ref{qaso}), 
q--Moyal algebras corresponding
to  (\ref{qsal}) and (\ref{qaal}) are obtained.
In \cite{odj} another 
q--star product is defined by using the coherent
states maps, however it does not lead to the algebras
which we deal with (\ref{qsal}),(\ref{qaal}).

Generalization of the Weyl
ordering (\ref{hsm}) to the q--phase space is not
obvious: a term of 
a monomial can be generalized 
by assuming that it is weighted with a factor  $q^\gamma ,$
where $\gamma$ is a number depending on the term
under consideration.
To emphasize the difficulties related to this ordering
procedure, let us suppose that there exist operators
$T_{m,n}$ leading to the ordinary Weyl ordered
operators in the $q=1$ limit,  satisfying
\be
[T_{m,n},T_{k,l}]_q =\sum_{r,s=0}^{A,B}
\Cc^{rs}_{mnkl}(\hbar ,q)
T_{r,s} ,
\ee
where for obtaining the correct classical limit
$\Cc$ should  satisfy
\be
\Cc^{(K)rs}_{mnkl}(0 ,q)=0.
\ee

An operator algebra is proposed in 
\cite{gf} as a generalization of (\ref{la1}) by replacing 
the factorial terms with the q--factorials:
\beqa
[\tbir_{m,n},\tbir_{k,l}]_q  & = &
\sum_{a=0}^b 
\sum_{c=0}^{2a+1} 
\fr{(-1)^c [m]! [n]! [k]! [l]!}{
[2a+1-b]! [b]! [m+c-2a-1]! [n-c]! [k-c]! [l+c-2a-1]! }  \nonumber \\
& & (i\hbar )^{2a+1}\tbir_{m+k-2a-1,n+l-2a-1} ,\lb{gfp}
\eeqa
where $b$ is given as in (\ref{la1}).
Symbol map independent of the definition of the 
q--Weyl ordered operators is
\be
\lb{gfs}
\Sc_W\tbir_{m,n} =p^mx^n,
\ee
so that, the $\star$--product reproducing
(\ref{gfp}), in terms of the q--Moyal bracket (\ref{qmoy}), is
\beqa
\star_W \equiv \sum_{M=0}^\infty 
(-\nu /2)^M\sum_{L=0}^M \fr{(-1)^L}{(M-L)! L!} 
(\dl_x x)^{M-L} 
(\dl_p p)^{L}  
\sum_{\al =0}^\infty (-i\hbar /2)^\al  \nonumber \\
\sum_{\beta =0}^\al
\fr{(-1)^\beta }{[\al -\beta ]![\beta ]!} 
\Dl_x^{\al -\beta } 
\Dl_p^{\beta }
\Dr_p^{\al -\beta }
\Dr_x^{\beta } 
(p \dr_p )^{M-L} 
(x \dr_x )^{L}  .
\eeqa
Obviously, independent of how the generalization is done,
we have
\be
T_{m,0}=P^m,\  
T_{0,m}=X^m.
\ee
Hence, by studying the q--commutators 
$[P^m,X^m]_q$ one can try to reach to
the other q--weighted monomials: e.g.
\be
\lb{e1}
P^2X^2-q^4X^2P^2 =i\hbar [2](PX +q^2XP),
\ee
suggests that $T_{1,1}\approx (PX +q^2XP). $
Similarly, explicit calculation gives
\be
P^2X^3-q^{6}X^3P^2= [2] (PX^2+q^4X^2P+q^2XPX),
\ee
yielding $T_{1,2}\approx (PX^2+q^4X^2P+q^2XPX).$
However, q--commutator of this operator with
$P:$
\beqa
P(PX^2+q^4X^2P+q^2XPX)-q^{2}(PX^2+q^4X^2P+q^2XPX)P
=  \nonumber \\
 i\hbar (1+q+q^2)
(PX+q^3XP),
\eeqa
suggests $T_{1,1}\approx (PX +q^3XP), $
which is in contradiction with the one suggested by (\ref{e1}).

In fact, having obstacles in defining q--Weyl
ordering is not surprising. If one considers monomials
reproducing the Weyl ordered ones in the $q=1$ limit
with a definite q weight they can not constitute a
complete basis. For having a complete basis one should
have monomials with all possible q weights.

\section{Applications:}

In terms of the q--star products and the related q--Moyal
brackets (\ref{qsso})--(\ref{qmoy}), 
one can study  quantum  as well as 
classical dynamics   on general grounds. 

If we deal with the $\star_S$--product,
q--classical dynamics can  be given in terms of the 
q--Poisson bracket defined for the observables 
$f(p,x)=\sum_if_i(p,x)$ and
$g(p,x)=\sum_jg_j(p,x)$ where $f_i$ and
$g_j$ are monomials in $p,x,$ as
\beqa
\{ f(x,p),g(x,p)\}_{q-P} &  \equiv  &  
\lim_{\hbar \rightarrow 0}\{ f(x,p),g(x,p)\}_{q-M}  \nonumber \\
& = & \sum_{i,j}[
q^{\al (f_i,g_j)}(D_p f_i) \exp (\nu \dl_ppx\dr_x)  (D_x g_j) \nonumber\\
&  & -q^{\al (g_j,f_i)}(D_p g_j) \exp (\nu \dl_ppx\dr_x)  (D_x f_i) ],
\lb{q-P}
\eeqa
where $\al (p^mx^n,p^kx^l)=nk.$ 
Thus, if $H=\sum_kH_k$ is the classical
hamiltonian where $H_k$ are monomials, 
equation of motion of the observable $f(p,x)$ is 
\be
\tau_q(f)=\{H,f\}_{q-P},
\ee
where $\lim_{q\rightarrow 1}\tau_q =d/dt.$
Now, observe that in general
\be
\{H,fg\}_{q-P}\neq f\{H,g\}_{q-P}+g\{H,f\}_{q-P},
\ee
which leads to
\be                           
\lb{ndf}
\tau_q(fg)\neq \tau (f)g+ f\tau (g).
\ee
Obviously, there are some exceptions like $f$ or $g$
is constant. One may think that $\star_S$
of $f$ and $g$ should be considered on the
left hand side of  (\ref{ndf}), however in
the limit 
$\hbar \rightarrow 0$ the q--star product will yield
$\star_S^{(q)}$ (\ref{nus}), and
$f\star_S^{(q)}g=q^{\al (f,g)}fg,$ so that there is not 
any difference.

In the ordinary 
classical mechanics canonical transformations leave
the basic Poisson brackets invariant. One of these 
is the point transformation  defined as
\be
\lb{kb}
u=f(x);\   p_u=\left(\del_xf(x)\right)^{-1}p  ,
\ee
where $f$ is an invertible function. For the 
q--classical mechanics
point transformation can be defined as
\be
\lb{qcan}
u=f(x);\   p_u=\left(D_xf(x)\right)^{-1}p  .
\ee
Now, in terms of the q--Poisson bracket (\ref{q-P}) one can
observe that
\be
\lb{trc}
\{u,p_u \}_{q-P} =-q^{\al (u,p_u)};\   \{p_u,u\}_{q-P}=1,
\ee
instead of the ones satisfied by $p,x:$
\be
\lb{xpp}
\{x,p\}_{q-P}=-q;\  \{p,x\}_{q-P}=1.
\ee
In (\ref{trc}) $\al (u,p_u)=x\del_x \log f(x)$ which is a number.
Let us have an example where 
$f(x)=\sqrt{x}:$
\be
\lb{neb}
u=\sqrt{x};\   p_u=\left[\fr{1}{2}\right]^{-1}\sqrt{x}p,
\ee
so that,
\be
\lb{trc1}
\{u,p_u \}_{q-P} =-q^{1/2};\   \{p_u,u\}_{q-P}=1.
\ee
In fact, the transformation 
(\ref{neb}) was studied in \cite{dd} and found that it
is a q--canonical transformation if the phase space operators
satisfy
\beqa
\hat{p}\hat{x}-q\hat{x}\hat{p} & = & i\hbar, \\
\hat{p}_u\hat{u}-\sqrt{q}\hat{u}\hat{p}_u & = & i\hbar,
\eeqa
which are consistent with (\ref{xpp}) and (\ref{trc1}).

When we deal with  q--quantum mechanics in the
Heisenberg picture,
time evolution  of an observable $f$ is  given by
\be
\lb{qem}
\tau (f)=\{H,f\}_{q-M}=\fr{1}{i\hbar }\sum_{i,j}\left(
q^{\al(H_i,f_j)}H_i\star f_j- q^{\al(f_j,H_i)}f_j\star H_i\right).
\ee                                                                   
Here $\star$ indicates one of the
q--star products (\ref{qsso})--(\ref{qaso}).

In the Schr\"{o}dinger picture time evolution of a 
time--dependent state vector $\psi$ is given by
\be
\lb{qse}
i\hbar \tau \left(\psi (t)\right)=\hat{H}\psi (t),
\ee
where $\hat{H}$ is the q--hamiltonian operator 
i.e. $\Sc (\hat{H})=H.$

In contrary to the ordinary quantum mechanics, in the 
q--deformed case relation between the Schr\"{o}dinger picture
(\ref{qse}) and the Heisenberg picture (\ref{qem})
is not clear. The unique common feature of the deformed
and non-deformed cases is the fact that in both of the 
cases symbols of the monomials are the same, namely $p^mx^n.$
The difference lies in the definition of the related star 
products. Thus, we may still assume that the symbol
of the evolution operator $\hat{U}(t)$ is
\be
U(t)=\Sc \left(\hat{U}(t)\right)=e^{\fr{it}{\hbar}H}.
\ee
Then we can adopt the definition of the path integral
of the ordinary time evolution 
given in terms of the star products and symbols in
\cite{ber} to define 
the q--path integral as
\be
G(t)=\lim_{N\rightarrow \infty}U(\fr{t}{N})
\star \cdots \star U(\fr{t}{N}).
\ee

When the canonical transformation (\ref{kb}) is performed
there will be some q--quantum corrections in the
transformed hamiltonian.
The kinetic term including the q--quantum corrections
can be studied in terms of the $\star_S$--product 
similar to the ordinary case\cite{du} 
by defining the 
transformed kinetic term as
\be
\lb{le}
\ti{H}_0=(D_xf)^{-1/2} \star_S p \star_S
(D_xf)^{-2}  \star_S D_xf \star_S p \star_S
(D_xf)^{-1/2}.
\ee
Because of the non-associativity of the q--star product
we should specify in which order  the multiplications will
be performed in (\ref{le}). However, we do not possess
a general procedure. 

As it is briefly illustrated, q--star products and q--Moyal
brackets are very useful in studying several aspects
of q--deformations on general grounds. However, the relations
to the other formulations of q--dynamics (e.g. see \cite{qdy})
and q--path integral definitions\cite{qpi} should be
studied.

\vspace{.5cm}

{\bf Acknowledgments}

A part of this work is done in II. Institut f\"{u}r 
Theoretische Physik, Uni Hamburg with
the support of A. von Humboldt Foundation. I thank them
for their kind hospitality.       
I am also thankful to M. Reuter for the discussions 
which were crucial for this work.


\vspace{2cm}

\newcommand{\jpa}{ J. Phys. A: Math. Gen. }
\newcommand{\prd}{ Phys. Rev. D } 
\newcommand{\plb}{ Phys. Lett. B }
\newcommand{\cmp}{ Commun. Math. Phys. }
\newcommand{\ijmp}{ Int. J. Mod. Phys A }
\newcommand{\jmp}{ J. Math. Phys. }
\newcommand{\mpl}{ Mod. Phys. Let. A }
\newcommand{\pr}{ Preprint }
\newcommand{\npb}{ Nucl. Phys. B }
\newcommand{\bi}{\bibitem}

\end{document}